\documentclass[8pt,twoside]{article}
\usepackage{eurosym}
\usepackage[utf8]{inputenc}
\usepackage{bbm}
\usepackage[english]{babel}
\usepackage[T1]{fontenc}
\usepackage{amsmath}
\usepackage{amsfonts}
\usepackage{amssymb}
\usepackage{amsthm}
\usepackage{anysize}
\usepackage{commath}

\usepackage{float}
\usepackage{array}
\usepackage{multirow}
\usepackage{tkz-tab}
\usepackage[colorlinks,linkcolor=black,citecolor=black]{hyperref}
\usepackage[cyr]{aeguill}
\usepackage{fancyhdr}
\usepackage{algpseudocode}
\usepackage{algorithm}

\addto\captionsfrench{}

\newtheorem{lemma}{Lemma}
\newtheorem{theorem}{Theorem}
\newtheorem{remark}{Remark}

\marginsize{40mm}{40mm}{20mm}{23mm}
\begin{document}

\title{\textbf{\ Primal-dual interior-point algorithm for linearly constrained convex optimization based  on a parametric  algebraic transformation }}
\date{}
\author{Aicha Kraria\\
Laboratory of Fundamental and Numerical Mathematics,\\
Department of Mathematics,\\
Setif 1 University-Ferhat Abbas, Algeria\\
aicha.kraria@univ-setif.dz
\and Bachir Merikhi\\
Laboratory of Fundamental and Numerical Mathematics,\\
Department of Mathematics,\\
Setif 1 University-Ferhat Abbas, Algeria
\and Djamel Benterki\\
Laboratory of Fundamental and Numerical Mathematics,\\
Department of Mathematics,\\
Setif 1 University-Ferhat Abbas, Algeria}
\maketitle

\pagestyle{fancyplain} 
\renewcommand{\sectionmark}[1]{\markright{\thesection:\ #1}{}}

\renewcommand{\baselinestretch}{1}\small \normalsize

\begin{abstract}
	In this paper, we present an interior point algorithm with a full-Newton step for solving a linearly constrained convex optimization problem, in which we propose a generalization of  the work of Kheirfam and Nasrollahi \cite{kheirfam2018full}, that consists in determining the descent directions through a parametric algebraic transformation.\\
	The work concludes with a complete study of the convergence of the algorithm and its complexity, where we show that the obtained algorithm achieves a polynomial complexity bounds.
	
\end{abstract}

\section{Introduction}

Interior point methods (IPMs) are among the most widely used and efficient methods for solving optimization problems. They first appeared in the 1950s, but developed and spread after the publication of  Karmarkar's paper in 1984 \cite{karmarkar1984new}, in which he proposed an interior point algorithm of  potential reduction type with polynomial complexity and very efficient for linear programming (LP). In the beginning of 1990's, a new class of interior point methods called the central path method appeared. This method of primal-dual type has a good theoretical and numerical behaviour, which has encouraged researchers to publish multiple studies on this method, in linear programming, convex quadratic programming, semidefinite programming and even for nonlinear programming.\\
To determine the descent directions of interior point methods, there are two main methods. One is based on kernel functions and requires two types of iterations, inner iterations and outer iterations (e.g. \cite{doi:10.1137/S1052623403423114, LI2015471, bouafia2016complexity, boudjellal2022complexity}). The other method was proposed by Darvay \cite{zs2002new} and is based on an equivalent algebraic transformation on the equation that defines the central path. This method has become widely used and so far a lot of research has been carried out on it. In 2003, Darvay \cite{darvay2003new} adopted this method where he proposed an interior point algorithm for (LP) based on the function $\psi(t)=\sqrt{t}$, and conducted a theoretical study on the convergence and complexity of the proposed algorithm. Subsequently, many researchers extended this work to various mathematical programs. The first extension from (LP) to convex quadratic programming was carried out in 2006 by Achache \cite{achache2006new}. Then, in 2008, Zhang et al. \cite{zhang2008new} extended this method to linearly constrained convex optimization (LCCO), where he showed that his algorithm had a polynomial complexity, namely $O(\sqrt{n}\log\frac{n}{\epsilon})$.\\
In their papers,  Wang and Bai \cite{wang2009primal}, Bai et al. \cite{bai2010polynomial}, Wang \cite{wang2012new}, Mansouri and Pirhaji \cite{mansouri2013polynomial} and  Asadi and Mansouri \cite{asadi2013polynomial} presented interior point algorithms based on the equivalent algebraic transformation technique to solve second-order cone optimization problem, convex quadratic semidefinite optimization problem, monotone linear complementarity problem over symmetric cones, monotone linear complementarity problems and $P_{*}(\kappa)$ horizontal linear complementarity problems, respectively, using the function $\psi(t)=\sqrt{t}$.\\
In 2016, Darvay et al. \cite{darvay2016complexity} proposed a new interior point algorithm in which they applied a new function $\psi(t)=t-\sqrt{t}$ in the equivalent algebraic transformation technique. Later, Kheirfam \cite{kheirfam2017new} and Darvay and Rig{\'o} \cite{darvay2018new} generalized this work to semidefinite optimization and symmetric optimization, respectively. Furthermore, several other interior point algorithms based mainly on Darvay's technique have been proposed using different functions, (see e.g. \cite{mansouri2012path, roumili2012weighted, kheirfam2013new, Mohamed2014AWF, asadi2017infeasible,  darvay2017, haddou2019generalized, billel2023interior, darvay2023interior, Zaoui2024}).\\
Recently, Kheirfam and Nasrollahi \cite{kheirfam2018full} enriched the analysis given by Darvay \cite{darvay2003new}, using the integer powers of the square root function. In this paper, we extend the work of Kheirfam and Nasrollahi \cite{kheirfam2018full} to (LCCO), in which we present lemmas that summarize our studies of convergence and complexity of the proposed algorithm, where we search for the integer that gives the best algorithmic complexity.\\
This paper is divided into five sections: in section \ref{s2}, we present the central path method for (LCCO). In section \ref{s3}, we extend Darvay's technique for determining descent directions to (LCCO), and we give the corresponding primal-dual algorithm. In section \ref{s4}, through some obtained results, we show the strict feasibility of all iterations, the convergence and the polynomial algorithmic complexity of the proposed algorithm. Finally, in section \ref{s5} we give a general conclusion.

\section{Central path method for (LCCO)}\label{s2}

In this work, we are interested in the following convex optimization problem:
\begin{equation*}\label{P}
\left\{ 
\begin{array}{c}
\min f(x), \\ 
Ax=b,\tag{P} \\ 
x\geq 0,%
\end{array}%
\right.
\end{equation*}
and its dual: 
\begin{equation*}\label{D}
\begin{cases}
\max b^{T}y+f(x)-x^{T}\nabla f(x),\notag \\ 
A^{T}y+z=\nabla f(x),\tag{D} \\ 
z\geq 0,\,y\in \mathit{R}^{m},\notag
\end{cases}
\end{equation*}
where $f:\,\mathit{R}^{n}\,\longrightarrow \mathit{R}$ is a convex and twice differentiable function, $A\in \mathit{R}^{m\times n}$ is a full rank matrix ($rank(A)=m<n$) and $b\in \mathit{R}^{m}$.\\
In the following, we assume that the sets of strictly feasible solutions of (\ref{P}) and (\ref{D}) are non-empty, i.e., there are $(\bar{x},\bar{y},\bar{%
z})$, such that: 
\begin{align*}
A\bar{x}& =b,\quad \bar{x}>0, \\
A^{T}\bar{y}+\bar{z}& =\nabla f(\bar{x}),\quad \bar{z}>0.
\end{align*}
Since $f$ is convex then, solving both problems (\ref{P}) and (\ref{D}) amounts to solving the following system, which is obtained by the optimality
conditions: 
\begin{align}
Ax& =b,\quad x\geq 0,  \notag  \label{11} \\
A^{T}y+z& =\nabla f(x),\quad z\geq 0, \\
xz& =0,  \notag
\end{align}
where $xz=(x_{1}z_{1},x_{2}z_{2},...,x_{n}z_{n})^{T}$.\\
To solve this system, we use the central path method, which consists in replacing the last equation of system (\ref{11}) by the equation $xz=\mu e_n$, where $\mu $ is a positive number and $e_{n}$ is a vector of $\mathit{R}^{n}$ in which all its elements are $1$. Thus, we obtain: 
\begin{align}
Ax& =b,\quad x>0,  \notag  \label{12} \\
A^{T}y+z& =\nabla f(x),\quad z>0, \\
xz& =\mu e_{n}.  \notag
\end{align}%
System (\ref{12}) has a unique solution $(x(\mu ),y(\mu ),z(\mu ))$ for each  $\mu >0$, this solution is called the $\mu $-center \cite{sonnevend71analytic}. The set of all $\mu $-center forms the central path of (\ref{P}) and (\ref{D}).\\
If $\mu \rightarrow 0$, then $\lim\limits_{\mu \rightarrow 0}(x(\mu ),y(\mu),z(\mu ))=(x^{\ast },y^{\ast },z^{\ast })$, i.e., the central path
converges to the optimal solutions of (\ref{P}) and (\ref{D}).\\
Finally, system (\ref{12}) can be written as follows: 
\begin{align}
Ax& =b,\quad x>0,  \notag  \label{13} \\
A^{T}y+z& =\nabla f(x),\quad z>0, \\
w^{2}& =e_{n},  \notag
\end{align}
where $w=\sqrt{\frac{xz}{\mu }}$.

\section{Darvay's technique for determining descent directions}\label{s3}

In this section, we use Darvay's technique \cite{darvay2003new} to determine the descent directions, where we replace the third equation in system (\ref{13}) by the equation $\psi (w^{2})=\psi (e_{n})$, where $\psi :\;]0,+\infty \lbrack \,\longrightarrow ]0,+\infty \lbrack $ is an invertible differentiable function, so we obtain: 
\begin{align}
Ax& =b,\quad x>0,  \notag  \label{14} \\
A^{T}y+z& =\nabla f(x),\quad z>0, \\
\psi (w^{2})& =\psi (e_{n}).  \notag
\end{align}
By applying Newton's method, we obtain the Newton's directions, which are solutions of system (\ref{15}): 
\begin{align}
A\Delta x& =0,  \notag  \label{15} \\
A^{T}\Delta y+\Delta z& =\nabla ^{2}f(x)\Delta x, \\
z\Delta x+x\Delta z& =\dfrac{\mu \bigl(\psi \bigl(e_{n})-\psi (w^{2})\bigr)}{%
\psi ^{^{\prime }}(w^{2})},  \notag
\end{align}%
and the new iteration of full Newton step is given by: 
\begin{equation*}
(x_{+},y_{+},z_{+})=(x,y,z)+(\Delta x,\Delta y,\Delta z).
\end{equation*}%
Now, we introduce the directions: 
\begin{equation*}
d_{x}=\frac{w\Delta x}{x},\,d_{z}=\frac{w\Delta z}{z},
\end{equation*}
then, system (\ref{15}) becomes: 
\begin{align}
\bar{A}d_{x}& =0,  \notag  \label{16} \\
\bar{A}^{T}\Delta _{y}+d_{z}& =Bd_{x}, \\
d_{x}+d_{z}& =p_{w},  \notag
\end{align}%
where 
\begin{align}
\bar{A}=\frac{1}{\mu }A\Bigl[diag\Bigl(\frac{x}{w}\Bigr)\Bigr]& ,\,B=\frac{1}{\mu }\Bigl[diag\Bigl(\frac{x}{w}\Bigr)\Bigr]\nabla ^{2}f(x)\Bigl[diag\Bigl(%
\frac{x}{w}\Bigr)\Bigr],  \notag  \label{17} \\
& p_{w}=\dfrac{\psi (e_{n})-\psi (w^{2})}{w\psi ^{^{\prime }}(w^{2})},
\end{align}%
such that $\frac{x}{w}=\bigl(\frac{x_{1}}{w_{1}},\frac{x_{2}}{w_{2}},...,\frac{x_{n}}{w_{n}}\bigr)^{T}$.

In the following, and inspired by the work of Kheirfam and Nasrollahi \cite{kheirfam2018full}, we propose the function $\psi (t)=t^{\frac{r}{2}}$, where $t>0$ and $r\in \mathit{N}^{\ast }$ for solving (LCCO) problem, which gives:
\begin{equation}
p_{w}=\dfrac{2e_{n}-2w^{r}}{rw^{r-1}},\,w>0.  \label{18}
\end{equation}
From this vector, we introduce a proximity measure $\Gamma :\mathit{R}_{++}^{n}\longrightarrow \,\mathit{R}_{+}$ defined by: 
\begin{equation}
\Gamma (x,z,\mu )=\Gamma (w)=\frac{\|p_w\|}{2}=\frac{1}{r}\biggl\|\dfrac{e_n-w^{r}}{w^{r-1}}\biggr\|,  \label{19}
\end{equation}
where $\mathit{R}_{++}^{n}=\{x\in \mathit{R}^{n}:x>0\}$ and $\|,\|$ denotes the Euclidean norm. 

\subsection{Primal-dual interior-point algorithm for (LCCO)}

\begin{algorithm}[H]
\caption{Primal-dual algorithm for (LCCO)}
\begin{algorithmic}
\State
\State \textbf{Input}
\State \quad  An accuracy parameter $\epsilon>0$;
\State \quad  An update parameter $0<\theta<1,\,\Bigl(\theta=\dfrac{1}{e^{2r}\sqrt{n}},\,r\in\mathbb{N}^*\Bigr)$;
\State \quad  A barrier parameter $\mu ^{(0)}=\frac{(x^{(0)})^Tz^{(0)}}{n}$;
\State \quad  An initial point  $(x^{(0)},y^{(0)},z^{(0)})$, such that $\Gamma (x^{(0)},z^{(0)},\mu ^{(0)})<\gamma=\dfrac{1}{e^r}$;
\State \textbf{begin}
\State \quad  $x:=x^{(0)}$; $y:=y^{(0)}$; $z:=z^{(0)}$; $\mu:=\mu ^{(0)}$; $k:=0;$
\State \quad \textbf{While} $x^Tz>\epsilon$ \textbf{do}
\State \qquad $\mu :=(1-\theta)\mu $;
\State \qquad Calculate $(\Delta x,\Delta y,\Delta z)$ from systems (\ref{15}) and (\ref{16}), then put:
\State \qquad $x:=x+\Delta x$;
\State \qquad $ y:=y+\Delta y$;
\State \qquad $ z:=z+\Delta z$;  $k:=k+1$;
\State \quad \textbf{end;}
\State \textbf{end.}
\end{algorithmic}
\label{algo}
\end{algorithm}
In order to facilitate the analysis of algorithm \ref{algo}, we introduce the following notation: $q_w=d_x-d_z,$ then: 
\begin{align}  \label{110}
\|p_w\|^{2}=\|d_x\|^{2}+\|d_z\|^{2}+2d_x^Td_z=\|q_w\|^{2}+4 d_x^{T}d_z.
\end{align}

\begin{remark}
Let $(d_x,\Delta y,d_z)$ be the solution of system (\ref{16}). We have: 
\begin{align}  \label{111}
d_x^{T}d_z=d_x^{T}Bd_x\ge 0,
\end{align}
because $B$ is a positive semidefinite matrix. 
\end{remark}
From (\ref{110}) and (\ref{111}), we deduce that: 
\begin{align}  \label{112}
\|p_w\|\ge \|q_w\|.
\end{align}

\section{Theoretical results}\label{s4}

In this section, we present our analysis of convergence and complexity of algorithm \ref{algo}, using the function $\psi(t)=t^\frac{r}{2},\,r\in\mathit{N}^*$.
\begin{lemma}
If $\Gamma(x,z,\mu)<1$ then the new iterations $x_+$ and $z_+$ are strictly feasible. 
\end{lemma}
\begin{proof}
Let $x_{+}(\alpha ),z_{+}(\alpha )$ be two vectors defined by: 
\begin{equation*}
x_{+}(\alpha )=x+\alpha \Delta x,\;z_{+}(\alpha )=z+\alpha \Delta z,
\end{equation*}%
such that $0\leq \alpha \leq 1$. We have: 
\begin{align}
\dfrac{1}{\mu}x_+(\alpha)z_+(\alpha)& =\frac{xz}{\mu}+\alpha\frac{(x\Delta z+z\Delta x)}{\mu}+\alpha^2\frac{\Delta x\Delta z}{\mu},\notag\\
&=w^{2}+\alpha w(d_x+d_z)+\alpha^2d_xd_z,\notag\\
&=w^{2}+\alpha w(d_x+d_z)+\alpha^2\biggl(\dfrac{p_w^{2}-q_w^{2}}{4}\biggr).\notag
\intertext{From (\ref{16}), we have:}
\dfrac{1}{\mu}x_+(\alpha)z_+(\alpha)&=(1-\alpha)w^2+\alpha(w^2+wp_w)+\alpha^2\biggl(\frac{p_w^{2}}{4}-\frac{q_w^{2}}{4}\biggr).\label{113}
\end{align}
Using (\ref{18}), we obtain: 
\begin{equation}\label{114}
w^{2}+wp_{w}=e_{n}+\dfrac{(r-2)w^{r+1}+2w-rw^{r-1}}{rw^{r-1}}.
\end{equation}%
As $w>0$, then after studying the function $\dfrac{(r-2)w^{r+1}+2w-rw^{r-1}}{rw^{r-1}}$ for each component of the vector $w$, we show easily that for all $r\in \mathit{N}^{\ast }$: 
\begin{equation}
w^{2}+wp_{w}\geq e_{n}-\frac{p_{w}^{2}}{4},  \label{115}
\end{equation}%
 this gives: 
\begin{align}
\dfrac{1}{\mu }x_{+}(\alpha )z_{+}(\alpha )& \geq (1-\alpha )w^{2}+\alpha \biggl(e_{n}-\frac{p_{w}^{2}}{4}\biggr)+\alpha ^{2}\biggl(\frac{p_{w}^{2}}{4}-\frac{q_{w}^{2}}{4}\biggr)\notag \\
& =(1-\alpha )w^{2}+\alpha \Biggl[e_{n}-\biggl((1-\alpha )\frac{p_{w}^{2}}{4}+\alpha \frac{q_{w}^{2}}{4}\biggr)\Biggr].\label{15a}
\end{align}%
Now, we will prove that $e_{n}-\biggl((1-\alpha )\dfrac{p_{w}^{2}}{4}+\alpha \dfrac{q_{w}^{2}}{4}\biggr)>0$. In fact, we have: 
\begin{align*}
\Biggl\|(1-\alpha)\frac{p_w^{2}}{4}+\alpha \frac{q_w^{2}}{4}\Biggr\|_\infty &\le (1-\alpha)\frac{\|p_w^{2}\|_\infty}{4}+\alpha \frac{\|q_w^{2}\|_\infty}{4},\\
&\le (1-\alpha)\frac{\|p_w\|^2}{4}+\alpha \frac{\|q_w\|^2}{4},\\
\intertext{and from (\ref{112}), we obtain:}
\Biggl\|(1-\alpha)\frac{p_w^{2}}{4}+\alpha \frac{q_w^{2}}{4}\Biggl\|_\infty&\le (1-\alpha)\frac{\|p_w\|^2}{4}+\alpha \frac{\|p_w\|^2}{4},\\
&= \frac{\|p_w\|^2}{4}=\Gamma^2.
\end{align*}
Consequently, if $\Gamma (x,z,\mu )<1$, then $\bigl\|(1-\alpha )\frac{p_{w}^{2}}{4}+\alpha \frac{q_{w}^{2}}{4}\bigr\|_{\infty }<1$, which means that $e_{n}-\biggl((1-\alpha )\dfrac{p_{w}^{2}}{4}+\alpha \dfrac{q_{w}^{2}}{4}\biggr)>0$, then from (\ref{15a}), we conclude that: 
\begin{equation*}
x_{+}(\alpha )z_{+}(\alpha )>0.
\end{equation*}
As $x_{+}(\alpha ),z_{+}(\alpha )$ are continuous functions, so they do not change sign on the interval $[0,1]$ and since $x_{+}(0)=x>0,\,z_{+}(0)=z>0$ then this implies that $x_{+}(1)=x+\Delta x=x_{+}>0,\,z_{+}(1)=z+\Delta z=z_{+}>0$, which proves the lemma.
\end{proof}

\begin{lemma}\label{lem2}
Let $w>0$ and $w_{+}=\sqrt{\dfrac{x_{+}z_{+}}{\mu }}$, then: 
\begin{equation*}
w_{+}\geq \sqrt{1-\Gamma ^{2}}e_{n},
\end{equation*}
where $\Gamma $ is defined in (\ref{19}).
\end{lemma}

\begin{proof}
When we take $\alpha=1$ in (\ref{113}), we obtain: 
\begin{align*}
w_+^2&=\dfrac{x_+z_+}{\mu}=w^2+wp_w+\frac{p_w^{2}}{4}-\frac{q_w^{2}}{4},
\intertext{and using (\ref{115}), we find:}
w_+^2&\ge e_n-\frac{q_w^{2}}{4},
\intertext{this implies that:}
(w_+)_i^{2}&\ge 1-\dfrac{\|q_w\|_{\infty}^{2}}{4},\quad \forall i=1:n,\\
&\ge 1-\dfrac{\|q_w\|^2}{4},\notag\\
\intertext{and from (\ref{112}), we deduce:}
(w_+)_i^{2}&\ge 1-\dfrac{\|p_w\|^2}{4}=1-\Gamma^2,\quad \forall i=1:n.
\end{align*}
Finally, we conclude that: 
\begin{equation*}
w_{+}\geq \sqrt{1-\Gamma ^{2}}e_{n}.
\end{equation*}
\end{proof}

\begin{lemma}\label{lem3}\cite{darvay2016complexity}\newline
Let's consider the decreasing function $f_{2}:]d,+\infty \lbrack \rightarrow \mathit{R}_{+}^{\ast }$, such that $d>0$. If we take the vector $w\in  \mathit{R}_{+}^{n}$, where $\min (w)=\min\limits_{i} (w_{i})\ge d$, then: 
\begin{equation*}
\Vert f_{2}(w)(e_{n}-w^{2})\Vert \leq f_{2}(\min (w))\Vert e_{n}-w^{2}\Vert\leq f_{2}(d)\Vert e_{n}-w^{2}\Vert .
\end{equation*}
\end{lemma}

\begin{lemma}\label{lem4}
If $\Gamma=\Gamma(x,z,\mu)<\dfrac{1}{e^r}$, then: 
\begin{align*}
\Gamma(x_+,z_+,\mu)\le \dfrac{e^r\Bigl(e^{r^2}-(e^{2r}-1)^{\frac{r}{2}}\Bigr)\Bigl((r-1)^2+1\Bigr)}{r(e^{2r}-1)^{\frac{r-1}{2}}}\Gamma^2,
\end{align*}
this means that the proximity measure decreases quadratically.
\end{lemma}

\begin{proof}
From (\ref{19}), we have: 
\begin{equation*}
\Gamma (x_{+},z_{+},\mu )=\frac{1}{r}\biggl\|\dfrac{\bigl(e_n-w_+^r\bigr)}{w_+^{r-1}\bigl(e_n-w_+^2\bigr)}(e_n-w_+^2)\biggr\|.
\end{equation*}
At the beginning, we calculate the expression $(e_{n}-w_{+}^{2})$ using (\ref{113}). 
\begin{align}
e_{n}-w_{+}^{2}& =e_{n}-\frac{x_{+}z_{+}}{\mu },  \notag \\
&=e_n-\biggr((w^2+wp_w)+\frac{p_w^2}{4}-\frac{q_w^2}{4}\biggl),\notag
\intertext{from (\ref{114}), we have:}
e_{n}-w_{+}^{2}& =e_{n}-\biggr(e_{n}+\dfrac{(r-2)w^{r+1}+2w-rw^{r-1}}{rw^{r-1}}+\frac{p_{w}^{2}}{4}-\frac{q_{w}^{2}}{4}\biggl),  \notag \\
& =\frac{q_{w}^{2}}{4}-\frac{p_{w}^{2}}{4}\biggr(\dfrac{4\bigl((r-2)w^{r+1}+2w-rw^{r-1}\bigr)}{rw^{r-1}p_{w}^{2}}+e_{n}\biggl),  \notag \\
& =\frac{q_{w}^{2}}{4}-\frac{p_{w}^{2}}{4}\biggr(\dfrac{(r-1)^{2}w^{2r}+(2r-2)w^{r}-r^{2}w^{2r-2}+e_{n}}{(e_{n}-w^{r})^{2}}\biggl).\label{116}
\end{align}
We can easily prove that: 
\begin{equation}
0\leq \dfrac{(r-1)^{2}w^{2r}+(2r-2)w^{r}-r^{2}w^{2r-2}+e_{n}}{(e_{n}-w^{r})^{2}}\leq (r-1)^{2}e_{n},  \label{117}
\end{equation}
for all $r\in \mathit{N}^{\ast },w>0\mbox{ and }w\neq e_{n}$.\newline
Using (\ref{116}) and (\ref{117}), we find: 
\begin{align}
\|{e_n-w_+^2}\|& =\biggl\|\frac{q_w^2}{4}-\frac{p_w^2}{4}\biggr(\dfrac{(r-1)^2w^{2r}+(2r-2)w^r-r^2w^{2r-2}+e_n}{(e_n-w^r)^2}\biggl)\biggr\|,  \notag \\
& \leq \frac{\|q_w^2\|}{4}+\biggl\|\frac{p_w^2}{4}\biggr(\dfrac{(r-1)^2w^{2r}+(2r-2)w^r-r^2w^{2r-2}+e_n}{(e_n-w^r)^2}\biggl)\biggr\|,  \notag \\
&\leq \frac{\|q_w\|^{2}}{4}+(r-1)^{2}\frac{\|p_w\|^{2}}{4},\notag
\intertext{and from (\ref{112}), we can write:}
\|e_n-w_+^2\|& \leq \frac{\|p_w\|^{2}}{4}+(r-1)^{2}\frac{\|p_w\|^{2}}{4}.  \notag
\intertext{Finally, and by using (\ref{19}), we find:}
\|e_n-w_+^2\|& \leq \bigl((r-1)^{2}+1\bigr)\Gamma ^{2}.  \label{118}
\end{align}

Now, we take the function $g_1(t)=\dfrac{\bigl(1-t^r\bigr)}{t^{r-1}\bigl(1-t^2\bigr)},t>0,t\neq 1$.\newline
After studying the variations of the function $g_1$, we find that $g_1^{\prime }(t)\le 0,\forall t>0,t\neq 1,$ which means that $g_1$ is a decreasing function.\newline
Using this, Lemma \ref{lem2} and Lemma \ref{lem3}, we obtain: 
\begin{align*}
\Gamma(x_+,z_+,\mu)&=\frac{1}{r}\bigl\|g_1(w_+)(e_n-w_+^2)\bigr\|, \\
&\le \frac{g_1(\sqrt{1-\Gamma^2})}{r}\bigl\|(e_n-w_+^2)\bigr\|,
\end{align*}
and as $\Gamma<\dfrac{1}{e^r}$, so $\sqrt{1-\Gamma^2}>\sqrt{1-\dfrac{1}{e^{2r}}}$, then, using this and (\ref{118}), we obtain: 
\begin{align*}
\delta(x_+,z_+,\mu)&\le \frac{g_1\biggl(\sqrt{1-\dfrac{1}{e^{2r}}}\biggr)}{r}\bigl((r-1)^2+1\bigr)\Gamma^2,
\intertext{thus, we conclude:}
\Gamma(x_+,z_+,\mu)&\le \dfrac{e^r\Bigl(e^{r^2}-(e^{2r}-1)^{\frac{r}{2}}\Bigr)\Bigl((r-1)^2+1\Bigr)}{r(e^{2r}-1)^{\frac{r-1}{2}}}\Gamma^2.
\end{align*}
\end{proof}

\begin{lemma}\label{lem5} 
If $\Gamma(x,z,\mu)<\dfrac{1}{e^r}$, then the
duality gap verifies the following inequality for all $n\in \mathit{N}^*$: 
\begin{align*}
(x_+)^Tz_+\le \mu \biggl(n+\dfrac{(r-1)^2}{e^{2r}}\biggr).
\end{align*}

\end{lemma}

\begin{proof}
From (\ref{113}), we have: 
\begin{align}
(x_{+})^{T}z_{+}& =\mu \sum\limits_{i=1}^{n}(w_{+})_{i}^{2},  \notag
\label{119} \\
& =\mu \sum\limits_{i=1}^{n}\biggl((w^{2}+wp_{w})_{i}+\frac{(p_{w}^{2})_{i}}{4}-\frac{(q_{w}^{2})_{i}}{4}\biggr).
\end{align}%
Using (\ref{114}), we obtain: 
\begin{align*}
w^{2}+wp_{w}& =e_{n}+\dfrac{(r-2)w^{r+1}+2w-rw^{r-1}}{rw^{r-1}}, \\
&=e_n+\frac{p_w^{2}}{4}\biggl(\dfrac{r(r-2)w^{2r}+2r w^r-r^2 w^{2r-2}}{(e_n-w^r)^2}\biggr).
\intertext{It is easy to prove that: $\dfrac{r(r-2)w^{2r}+2r w^r-r^2 w^{2r-2}}{(e_n-w^r)^2}\le r(r-2)e_n$, this implies that:}
w^2+wp_w&\le e_n+r(r-2)\frac{p_w^{2}}{4}.
\end{align*}
Using this and (\ref{119}), we obtain: 
\begin{align*}
(x_{+})^{T}z_{+}& \leq \mu \sum\limits_{i=1}^{n}\biggl(1+r(r-2)\frac{(p_{w}^{2})_{i}}{4}+\frac{(p_{w}^{2})_{i}}{4}-\frac{(q_{w}^{2})_{i}}{4}%
\biggr), \\
& =\mu \biggl(n+(r-1)^{2}\frac{\|p_w\|^{2}}{4}-\frac{\|q_w\|^{2}}{4}\biggr), \\
&\le \mu \biggl(n+(r-1)^2\Gamma^2\biggr).
\intertext{Since $\Gamma<\dfrac{1}{e^r}$, then:}
(x_{+})^{T}z_{+}& \leq \mu \biggl(n+\dfrac{(r-1)^{2}}{e^{2r}}\biggr).
\end{align*}
\end{proof}

\begin{lemma}
If $\Gamma(x,z,\mu)<\gamma=\dfrac{1}{e^r}$ and $\theta=\dfrac{1}{e^{2r}\sqrt{n}}$, then: 
\begin{align*}
\Gamma(x_+,z_+,\mu_+)<\dfrac{1}{e^r},
\end{align*}
for all $r\in \mathit{N}^*$.
\end{lemma}

\begin{proof}
We assume that $\bar{w}=\sqrt{\dfrac{x_+z_+}{\mu_+}}$, then: 
\begin{align*}
\bar{w}=\sqrt{\dfrac{x_+z_+}{(1-\theta)\mu}}=\frac{1}{\sqrt{1-\theta}}w_+>0,\quad (0<\theta<1).
\end{align*}
From (\ref{19}),we have: 
\begin{align}  \label{120}
\Gamma(\bar{w})&=\frac{\|p_{\bar{w}}\|}{2},  \notag \\
&=\frac{1}{r}\biggl\|\dfrac{e_n-\bar{w}^{r}}{\bar{w}^{r-1}(e_n-\bar{w}^2)}(e_n-\bar{w}^2)\biggr\|, \\
&=\frac{1}{r}\bigl\|g_1(\bar{w})(e_n-\bar{w}^2)\bigr\|.  \notag
\end{align}
In Lemma \ref{lem4}, we showed that $g_1(t)$ is a decreasing function for all $t>0,t\neq 1$, so using this result and Lemma \ref{lem3},
we obtain: 
\begin{align*}
\Gamma(\bar{w})&\le\frac{1}{r}g_1(\min\bar{w})\bigl\|e_n-\bar{w}^2\bigr\|,
\intertext{and since}
\bar{w}&=\frac{1}{\sqrt{1-\theta}}w_+,
\intertext{then, from Lemma \ref{lem2}, we can write:}
\bar{w}&\ge \frac{\sqrt{1-\Gamma^2}}{\sqrt{1-\theta}}e_n,
\intertext{which is equivalent to:}
\min(\bar{w})&\ge \frac{\sqrt{1-\Gamma^2}}{\sqrt{1-\theta}},
\end{align*}
this implies that: 
\begin{align}
\Gamma(\bar{w})&\le\frac{1}{r}g_1\biggl(\frac{\sqrt{1-\Gamma^2}}{\sqrt{1-\theta}}\biggr)\bigl\|e_n-\frac{1}{1-\theta}w_+^2\bigr\|,\notag \\
\intertext{from (\ref{118}), we obtain:}
\Gamma(\bar{w})&\le \frac{1}{r(1-\theta)}g_1\biggl(\frac{\sqrt{1-\Gamma^2}}{\sqrt{1-\theta}}\biggr)\biggl(\theta\sqrt{n}+((r-1)^2+1)\Gamma^2\biggr).\label{121}
\end{align}
We will now discuss three cases, depending on the value of $r$, using the values $\gamma=\dfrac{1}{e^r}$ and $\theta=\dfrac{1}{e^{2r}\sqrt{n}}$. 

\subsubsection*{\textit{\protect Case 1: $r=1$}}

In this case we have $\Gamma<\dfrac{1}{e}$ and $\theta=\dfrac{1}{e^{2}\sqrt{n}}$, so 
\begin{equation*}
\begin{cases}
\sqrt{1-\Gamma^2}>\sqrt{1-\dfrac{1}{e^{2}}}=a_1, \\ 
\sqrt{1-\theta}>\sqrt{1-\dfrac{1}{e^{2}}}=a_1,\;\mbox{for all $n>1$},
\end{cases}
\end{equation*}
this gives: 
\begin{align*}
\Gamma(\bar{w})&\le \frac{1}{(1-\theta)}g_1\biggl(\frac{a_1}{\sqrt{1-\theta}}\biggr)\biggl(\theta\sqrt{n}+\Gamma^2\biggr), \\
&=\frac{1}{\sqrt{1-\theta}}\biggl(\dfrac{1}{\sqrt{1-\theta}+a_1}\biggr)\biggl(\theta\sqrt{n}+\Gamma^2\biggr), \\
&<\frac{1}{a_1}\biggl(\dfrac{1}{2a_1}\biggr)\biggl(\dfrac{1}{e^{2}}+\dfrac{1}{e^2}\biggr), \\
&<\dfrac{1}{e}.
\end{align*}
In the first case, we proved that if $\Gamma(x,z,\mu)<\dfrac{1}{e}$ and $\theta=\dfrac{1}{e^{2}\sqrt{n}}$, then $\Gamma(x_+,z_+,\mu_+)=\Gamma(\bar{w})<\dfrac{1}{e}$.

\subsubsection*{\textit{\protect Case 2: $r=2$}}

We have: $\Gamma<\dfrac{1}{e^2}$ and $\theta=\dfrac{1}{e^{4}\sqrt{n}}$, so 
\begin{equation*}
\begin{cases}
\sqrt{1-\Gamma^2}>\sqrt{1-\dfrac{1}{e^{4}}}=a_2, \\ 
\sqrt{1-\theta}>\sqrt{1-\dfrac{1}{e^{4}}}=a_2,\;\mbox{for all $n>1$},
\end{cases}
\end{equation*}
this yields: 
\begin{align*}
\Gamma(\bar{w})&\le \frac{1}{2(1-\theta)}g_1\biggl(\frac{a_2}{\sqrt{1-\theta}}\biggr)\biggl(\theta\sqrt{n}+2\Gamma^2\biggr), \\
&=\frac{1}{2\sqrt{1-\theta}}\biggl(\dfrac{1}{a_2}\biggr)\biggl(\theta\sqrt{n}+2\Gamma^2\biggr), \\
&<\frac{1}{2a_2}\biggl(\dfrac{1}{a_2}\biggr)\biggl(\dfrac{1}{e^{4}}+\dfrac{2}{e^4}\biggr), \\
&<\dfrac{1}{e^2}.
\end{align*}
Also in this case, we proved that if $\Gamma(x,z,\mu)<\dfrac{1}{e^2}$ and $\theta=\dfrac{1}{e^{4}\sqrt{n}}$, then $\Gamma(x_+,z_+,\mu_+)=\Gamma(\bar{w})<\dfrac{1}{e^2}$.

\subsubsection*{\textit{\protect Case 3: $r\ge 3$}}

 Now, we take $\Gamma(x,z,\mu)< \dfrac{1}{e^r}$ and $\theta=\dfrac{1}{e^{2r}\sqrt{n}}$, such that $r\ge 3$, so we obtain: 
\begin{align}
\begin{cases}
\sqrt{1-\Gamma^2}>\sqrt{1-\dfrac{1}{e^{2r}}}\ge \sqrt{1-\dfrac{1}{e^{6}}}=a_3, \\ 
\sqrt{1-\theta}>\sqrt{1-\dfrac{1}{e^{2r}}}=\dfrac{\sqrt{e^{2r}-1}}{e^{r}},\;\mbox{for all $n>1$},\label{122}
\end{cases}
\end{align}
because $\sqrt{1-\dfrac{1}{e^{2r}}}$ is an increasing function for all $r\ge 3$. From this and (\ref{121}) we can write: 
\begin{align*}
\Gamma(\bar{w})&\le \frac{1}{r(1-\theta)}g_1\biggl(\frac{a_3}{\sqrt{1-\theta}}\biggr)\biggl(((r-1)^2+1)\Gamma^2+\theta\sqrt{n}\biggr), \\
&<\frac{1}{3\sqrt{1-\theta}}\biggl(\dfrac{(1-\theta)^{\frac{r}{2}}-a_3^r}{a_3^{r-1}((1-\theta)-a_3^2)}\biggr)\biggl(\dfrac{(r-1)^2+1}{e^{2r}}+\dfrac{1}{e^{2r}}\biggr).
\end{align*}
The function $g_2(\sqrt{1-\theta})=\dfrac{(1-\theta)^{\frac{r}{2}}-a_3^r}{a_3^{r-1}((1-\theta)-a_3^2)}$ is an increasing function for all $\sqrt{1-\theta}<1$. This results in: 
\begin{align*}
\Gamma(\bar{w})&<\frac{1}{3\sqrt{1-\theta}}\biggl(\dfrac{1-a_3^r}{a_3^{r-1}(1-a_3^2)}\biggr)\biggl(\dfrac{(r-1)^2+1}{e^{2r}}+\dfrac{1}{e^{2r}}\biggr),
\intertext{using (\ref{122}), we obtain:}
\Gamma(\bar{w})&<\frac{e^{r}}{3\sqrt{e^{2r}-1}}\biggl(\dfrac{1-a_3^r}{a_3^{r-1}(1-a_3^2)}\biggr)\biggl(\dfrac{(r-1)^2+1}{e^{2r}}+\dfrac{1}{e^{2r}}\biggr), \\
&=\frac{1}{e^{r}}\Biggl(\frac{e^{r}}{3\sqrt{e^{2r}-1}}\biggl(\dfrac{1-a_3^r}{a_3^{r-1}(1-a_3^2)}\biggr)\biggl(\dfrac{(r-1)^2+1}{e^{r}}+\dfrac{1}{e^{r}}\biggr)\Biggr).
\end{align*}
After some algebraic calculations, we find that: 
\begin{align*}
&\frac{e^{r}}{3\sqrt{e^{2r}-1}}<1
\intertext{and}
&\biggl(\dfrac{1-a_3^r}{a_3^{r-1}(1-a_3^2)}\biggr)\biggl(\dfrac{(r-1)^2+1}{e^{r}}+\dfrac{1}{e^{r}}\biggr)<1,
\end{align*}
which means that: 
\begin{align*}
\Gamma(\bar{w})&<\frac{1}{e^{r}}.
\end{align*}
Finally, the last case asserts that if $r\ge 3,\,\Gamma(x,z,\mu)<\dfrac{1}{e^r}$ and $\theta=\dfrac{1}{e^{2r}\sqrt{n}}$, then $\Gamma(x_+,z_+,\mu_+)<\dfrac{1}{e^r}$, this completes the proof.
\end{proof}

\begin{lemma}\label{lem7}
Let $(x^{(0)},y^{(0)},z^{(0)})$ be a strictly primal-dual feasible initial point, $\Gamma(x^{(0)},z^{(0)},\mu^{(0)})<\dfrac{1}{e^r},\,\epsilon>0$, $\mu^{(0)}=\frac{(x^{(0)})^Tz^{(0)}}{n}$, $0<\theta<1$ and $x^{(k+1)},z^{(k+1)}$ are the vectors obtained after $k+1$ iterations. If 
\begin{equation*}
k\ge \Biggl[\frac{1}{\theta}\log\dfrac{\mu^{(0)}\bigl(n+\frac{(r-1)^2}{e^{2r}}\bigr)}{\epsilon}\Biggr],
\end{equation*}
then $\bigl(x^{(k+1)}\bigr)^Tz^{(k+1)}\leq\epsilon $, for all $r\in\mathit{N}^*$.
\end{lemma}

\begin{proof}
From Lemma \ref{lem5}, we have: 
\begin{align*}
\bigl(x^{(k+1)}\bigr)^Tz^{(k+1)}&\le \mu^{(k)} \biggl(n+\dfrac{(r-1)^2}{e^{2r}}\biggr), \\
&=(1-\theta)^k \mu^{(0)}\biggl(n+\dfrac{(r-1)^2}{e^{2r}}\biggr),
\end{align*}
because $\mu^{(k)}=(1-\theta)^k \mu^{(0)}.$ This implies that to get
\begin{align*}
\bigl(x^{(k+1)}\bigr)^Tz^{(k+1)}\leq\epsilon
\intertext{ it is sufficient to verify that: }
(1-\theta)^k\mu^{(0)}\biggl(n+\dfrac{(r-1)^2}{e^{2r}}\biggr)\leq\epsilon.
\end{align*}
This last inequality is satisfied if 
\begin{align*}
&\log(1-\theta)^k+\log\mu^{(0)}\biggl(n+\dfrac{(r-1)^2}{e^{2r}}\biggr)\le
\log\epsilon, \\
\Longleftrightarrow & -k\log(1-\theta)\ge \log\dfrac{\mu^{(0)}\biggl(n+\dfrac{(r-1)^2}{e^{2r}}\biggr)}{\epsilon}.
\end{align*}
Moreover, we know that $-\log(1-\theta)\ge \theta,\,\forall\, 0<\theta< 1$,
so $-k\log(1-\theta)\ge k\theta,$ and therefore: 
\begin{align*}
-k\log(1-\theta)&\ge \log\dfrac{\mu^{(0)}\biggl(n+\dfrac{(r-1)^2}{e^{2r}}\biggr)}{\epsilon}, \\
\intertext{if}
k\theta&\ge \log\dfrac{\mu^{(0)}\biggl(n+\dfrac{(r-1)^2}{e^{2r}}\biggr)}{\epsilon}, \\
\intertext{hence}
k&\ge\frac{1}{\theta} \log\dfrac{\mu^{(0)}\biggl(n+\dfrac{(r-1)^2}{e^{2r}}\biggr)}{\epsilon}.
\end{align*}
\end{proof}

\begin{theorem}
If we take $\theta=\dfrac{1}{e^{2r}\sqrt{n}}$, the algorithm \ref{algo} converges to the optimal solution of (\ref{P}) and (\ref{D})
after 
\begin{align*}
\Biggl[ e^{2r}\sqrt{n}\log\dfrac{\mu^{(0)}\biggl(n+\dfrac{(r-1)^2}{e^{2r}}\biggr)}{\epsilon}\Biggr]
\end{align*}
iterations.
\end{theorem}

\begin{proof}
We replace $\theta$ by $\dfrac{1}{e^{2r}\sqrt{n}}$ in Lemma \ref{lem7}, we obtain the result.
\end{proof}

\section{General conclusion}\label{s5}
In this paper, we have successfully extended the work of Kheirfam and Nasrollahi \cite{kheirfam2018full} to (LCCO), by proposing a class of directions based on the parametric function $\psi(t)=t^{\frac{r}{2}}$, where $r\in \mathit{N}^*$. We have also demonstrated the convergence of the proposed algorithm and the theoretical polynomial complexity.  Finally, we found that the best complexity is given for $r=1$, i.e., when $\psi(t)=\sqrt{t}$.

\phantomsection
\bibliographystyle{abbrv}
\bibliography{Aicha}

\end{document}